\documentclass[12pt]{amsart}
\usepackage{xspace,amssymb,amsfonts,amsthm}
\usepackage{graphicx,euscript,eufrak}
\usepackage[all]{xy}
\title[Wall-Crossing Morphisms in Khovanov-Rozansky Homology]{Wall-Crossing Morphisms in\\
 Khovanov-Rozansky Homology}
\author{Nadya Shirokova}
\address{Department of Mathematics\\
Stanford University\\
Palo Alto, CA}
\author{Ben Webster}
\address{School of Mathematics\\
  Institute for Advanced Study\\
  Princeton, NJ 08540}
\email{bwebste@math.berkeley.edu}
\urladdr{http://math.berkeley.edu/~{}bwebste}
\subjclass[2000]{57M25}
\begin{document}
\begin{abstract}
  We define a wall-crossing morphism for Khovanov-Rozansky homology;
  that is, a map between the KR homology of knots related by a
  crossing change.  Using this map, we extend KR homology to an
  invariant of singular knots categorifying the Vasilliev derivative
  of the HOMFLY polynomial, and of $\mathfrak{sl}_n$ quantum
  invariants.
\end{abstract}
\maketitle
 \newtheorem{theorem}{Theorem}[section]
  \newtheorem*{theorem*}{Theorem}
  \newtheorem{proposition}[theorem]{Proposition}
  \newtheorem{lemma}[theorem]{Lemma}
  \newtheorem{corollary}[theorem]{Corollary}
  \newtheorem{conjecture}{Conjecture}

 \newtheorem{defi}{Definition}[section]
\newtheorem{defith}[theorem]{Definition/Theorem}

 \newcommand{\nc}{\newcommand}
 \newcommand{\renc}{\renewcommand}
 \nc{\df}{\bf}
 \nc{\ti}{\tilde}
  \nc{\al}{\alpha}
  \nc{\mcC}[2]{\mathcal{C}^{#1}_{#2}}
  \nc{\ep}{\epsilon}
  \nc{\maf}[2]{\mathsf{MF}^{#1}_{#2}}
  \nc{\MF}{\mathbf{MF}}
  \nc{\KR}{\mathcal{KR}}
  \nc{\KRh}{\mathcal{K}}
  \nc{\Kom}{\mathrm{Kom}}
  \nc{\mc}[1]{\mathcal{#1}}
  \nc{\Z}{\mathbb{Z}}
 \nc{\C}{\mathbb{C}}
  \nc{\vp}{\varphi}
  \nc{\By}{\mathbf{y}}
  \nc{\Bx}{\mathbf{x}}
  \nc{\Bt}{\mathbf{t}}
  \nc{\No}{N_1}
  \nc{\Nt}{N_2}
  \nc{\Kmf}[2]{\ti{\mc Z}_{#1,#2}}
  \nc{\mf}{matrix factorization\xspace}
  \nc{\mfs}{matrix factorizations\xspace}
  \nc{\qpm}{\textit{q-p}-morphism\xspace}
  \nc{\qpim}{\textit{q-p}-isomorphism\xspace}
  \nc{\qpms}{\textit{q-p}-morphisms\xspace}
  \nc{\dpl}{d_+}
  \nc{\dm}{d_-}
  \nc{\dpm}{d_\pm}
  \nc{\fr}[1]{\mathfrak{#1}}
  \nc{\Ht}{H^t}  
  \nc{\fH}{\EuScript{H}}
  \nc{\nai}{\EuScript{N}}
 \nc{\excise}[1]{}
 \nc{\Hhv}[2]{H_{hv}^{#1}\left(#2\right)}
 \nc{\Hv}[2]{H_{h}^{#1}\left(#2\right)}
 \nc{\hd}[1]{\EuScript{K}(#1)}
 \nc{\Span}{\mathrm{span}\,}
 \nc{\Hh}[2]{H_{v}^{#1}\left(#2\right)}
\nc{\KRH}{Khovanov-Rozansky homology\xspace}
\nc{\la}{\lambda}
\nc{\si}{\sigma}
\nc{\wc}{\EuScript{W}}
\nc{\doublepoint}{K_0}
\nc{\overcrossing}{K_+}
\nc{\undercrossing}{K_-}
\nc{\HE}{\EuScript{H}}
\nc{\A}{\mc A}
\nc{\KA}{\mc K(\A)}
\nc{\too}{\longrightarrow}
\nc{\ad}{\al}
\nc{\Q}{\mathbb{Q}}
\nc{\HH}{H\! H}
\nc{\Ext}{\mathrm{Ext}}
\nc{\VKR}{\mathcal{VKR}}

\section{Introduction}
\label{sec:introduction}
In \cite{Shi07a}, the first author outlined a program of
classification of combinatorially defined homological theories, such
as Khovanov and Ozsvath-Szabo's categorifications of knot and
3-manifold invariants.
 
In this paper we construct the wall-crossing morphism for
Khovanov-Rozansky (KR) homology which will allow us to put the
KR theory into such a framework.

The theory that was outlined in \cite{Shi07a} combines results of V. Vassiliev,
A. Hatcher and M. Khovanov and the resulting theory can be considered as
a ``categorification of Vassiliev theory'' or a classification of
categorifications of knot invariants. We defined Khovanov homology for singular knots, introduced the definition of a
theory of finite type $n$ and have shown that Khovanov homology restricted
to the subcategories of knots of bounded crossing number have finite
type \cite{Shi07b}.

The main idea in \cite{Shi07a} was to consider a knot homology theory as
a local system, or a constructible sheaf on the space of all objects
(knots, including singular ones), extend this local system to the
singular locus and introduce the analogue of the ``Vassiliev derivative''
for categorifications.
    
The Khovanov homology was just the first example of a theory satisfying
our axioms and in the present paper we show that the analogous
constructions can be carried out for Khovanov-Rozansky homology.
   

In \cite{Vas92}, Vassiliev introduced finite type invariants by considering the
space of all immersions of $S^1$ into $R^3$ and relating the topology of
the singular locus to the topology of its complement via Alexander
duality. He resolved and cooriented the discriminant (the space of
immersions with self-intersections) and introduced a spectral sequence
with a filtration, which suggested the simple geometrical and
combinatorial definition of an invariant of finite type.

Let $\lambda$ be an arbitrary invariant of oriented knots in oriented
space with values in an abelian group $A$. Extend $\lambda$ to be an
invariant of $1$-singular knots (knots that may have a single
singularity which is a double point), using the formula which was
interpreted by Birman and Lin as a ``Vassiliev derivative'': if $K_0$
has a single double point and $K_+$ and $K_-$ are its ``positive'' and
``negative'' resolutions, then 
\begin{equation}
  \label{eq:1}
\lambda(\doublepoint)=\lambda(\overcrossing)-\lambda(\undercrossing)  
\end{equation}
Furthermore, extend $\lambda$ to the set $\mathcal K^n$ of $n$-singular
knots (knots with $n$ double points) by repeatedly using the 
relation \eqref{eq:1}. 
\begin{defi}
We say that $\lambda$ is of type
$n$ if its extension $\left.\lambda\right|_{\mathcal K^{n+1}}$ to
$(n+1)$-singular knots vanishes identically. We say that $\lambda$ is of
{\bf finite type} if it is of type $n$ for some $n$.    
\end{defi}
Given this formula, the definition of an invariant of finite type n
becomes similar to that of a polynomial: its $(n+1)$-st Vassiliev
derivative is zero.
 
All known invariants are either of finite type, or are
infinite linear combinations of those.  For example, it was shown by
Bar-Natan that the nth coeffitient of the Conway polynomial is a
Vassiliev invariant of order $n$.
 
To categorify the Vassiliev derivative for a link homology theory
$\HE$, where $\HE(K)$ is a complex attached to a link $K$ (typically whose Euler characteristic is a known link invariant), we must
define a {\bf wall-crossing morphism} $\wc_w:\HE(K_-)\to \HE(K_+)$ for
each wall $w$ of codimension 1 in the discriminant.

We let $\HE(K_0)$ be the cone of $\wc_w$:
\begin{equation*}
  \HE(K_0)\cong \HE(K_+)\oplus \HE(K_-)[1]\hspace{10mm} d_0=
  \begin{bmatrix}
    -d_+&\wc_w\\0& -d_-
  \end{bmatrix}
\end{equation*}

Assuming these satisfy the obvious commutativity relation
$\wc_{w_1}\!\wc_{w_2}=\wc_{w_2}\!\wc_{w_1}$ for any two walls
$w_1,w_2$ which meet, then $\wc$ allows us to extend $\HE$ to the
whole discriminant, by simply taking successive cones, as described in Section~\ref{sec:the categorification of the Vassiliev derivative}.

Our principal result in this paper is to define just such a wall
crossing morphism for KR homology $\KR_N$ for any integer $N$, or HOMFLY homology, which we denote $\KR_\infty$.
\begin{theorem}
  For $N=1,\ldots, \infty$ and for any wall $w$ in the discriminant, we have a map
  $\wc_w:\KR_N(K_-)\to\KR_N(K_+)$.  The wall-crossing maps for adjacent
  walls commute, so this defines an extension of $\KR_N$ to the
  discriminant locus.
  
  The Euler characteristic of this extension is the extension of the
  $\mathfrak{sl}_N$-quantum invariant (or HOMFLY polynomial) to the
  discriminant by Vassiliev derivative.
\end{theorem}
There is, of course, an obvious guess for this wall-crossing map when
$N <\infty$: that defined on $\KR_N$ by the unique cobordism of
minimal genus between $K_+$ and $K_-$.  Unfortunately, the definition
of Khovanov and Rozansky \cite{KR1} is only well defined up to
constant factors.  Since it very important that our wall-crossing
morphisms commute (and not, say, anti-commute), we will need to
construct them explicitly, which is, in fact, much simpler than
attempting to compute the action of a cobordism.  Unfortunately, our
wall-crossing maps are also only defined up to scalar, but it is clear
from the construction that the scalars can be chosen consistently so
that wall-crossing maps commute.  This approach has the further
advantage of applying to the triply graded theory introduced in
\cite{KR2}, where even projective functoriality has not been
described.

It seems likely to the authors that this wall-crossing map has
connections to the action of the braid cobordisms on the derived
category of coherent sheaves on the cotangent bundle of the flag
variety described by Khovanov and Thomas \cite{KT}, but at the moment
it seems unclear how.

One family of algebraic invariants of singular knots has already been
categorified (as part of the definition of Khovanov-Rozansky
homology): the MOY invariants, which give one extension of the HOMFLY
polynomial and $\mathfrak{sl}_n$ quantum invariants to the
discriminant locus, one which is different from the Vassiliev
derivative.  While this categorification and our theory both give
extensions of \KRH to singular knots, they categorify different
extensions of the HOMFLY polynomial to singular knots, and thus
obviously differ.

\section*{Acknowledgments}
\label{sec:acknowledgments}

The authors would like to thank J\o rgen Ellegaard Andersen and
Nicolai Reshetikhin for their efforts toward making the Center for
Topology and Quantization of Moduli Spaces a great place to do
mathematics, and in particular, for organizing the retreat where this
collaboration began in November 2006.

B.W. was supported under a National Science Foundation Graduate Research
Fellowship, the National Science Foundation RTG grant DMS-0354321, and
a Danish National Research Foundation Niels Bohr Professorship grant.

\section{Khovanov-Rozansky Homology}
\label{sec:khov-rozansky-homol}

Since it is somewhat involved, we will only sketch the definition of
\KRH here, leaving out details that will not be relevant for our
argument. For a more complete definition, see the papers of
Rasmussen~\cite{Ras06}, Webster~\cite{Web06a}, Khovanov~\cite{Kho} or
the original papers of Khovanov and Rozansky~\cite{KR1}.  We will
concentrate on reduced HOMFLY homology, though our results are equally
valid for unreduced homology.

Reduced HOMFLY homology is defined as follows: 

Using Vogel's algorithm, write your knot $K$ as the closure of a braid
$\si=\si_{i_1}^{\ep_1} \cdots \si_{i_m}^{\ep_m}$, where $\{\si_i\}$
are the standard generators of the braid group $B_n$.  We let
$\ad(\si)=\frac 12\left(n-\sum_{j=1}^m\ep_{i_j}\right)$, that is,
half the braid index minus half the writhe.

Let $S=\Q[x_1,\ldots, x_n]/(x_1+\cdots+x_n)$, equipped with the
obvious $S_n$-action, and considered with the grading where
$\deg(x_i)=2$.  If $s_i=(i,i+1)$, then $S^{s_i}$ is the subring
generated by $x_j$ for $j\neq i,i+1$, and the symmetrized polynomials
$x_i+x_{i+1}$ and $x_ix_{i+1}$.

We define $S_{i}=S\otimes_{S^{s_i}}S\{-1\}$ (where $\{a\}$ denotes
degree shift by $a$).  This is free of rank two for the left and right
action of $S$, and is equipped with an obvious multiplication map
$m_i:S_i\to S$.  Less obviously, we have a map $\iota_i:S\{2\}\to
S_i$, such that $\iota_i(1)=x_i\otimes 1-1\otimes x_{i+1}$.

We can define a categorification $F$ of the braid group setting
\begin{align*}
  F(\si_i)&=\cdots\too 0 \too \,\,\,S\{2\}\,\,\, \too \,S_i\{1\}\,\, \too 0\too \cdots \\
  F(\si_i^{-1})&=\cdots\too 0 \too S_i\{-1\} \too S\{-2\} \too 0\too \cdots 
\end{align*}
and extending to an arbitrary element of $B_n$ by the relation $F(\si\si')=F(\si)\otimes_R F(\si')$.  By a theorem of Rouquier \cite{Rou04}, this is well-defined up to homotopy equivalence of complexes.

The triply-graded HOMFLY homology can be constructed by applying the functor
Hochschild homology $\HH^*$ to $F(\si)$.  For each $i$, we have a
complex $\EuScript{F}^i(\si)=\HH^i(F(\si))$ (with Hochschild homology applied termwise) of
graded modules over $S$.
\begin{defi}
  The {\bf HOMFLY homology} of $\bar{\si}$ is the doubly graded complex of $S$-modules
  whose graded pieces are
  $\KR^{i,j,k}(\bar{\si})=\EuScript{F}^{i-\ad(\si)}(\si)_{j}$, that is, the elements of 
  grade $j$ of the complex
  $\HH^{i-\ad(\si)}(F(\si))$.  
\end{defi}
As Khovanov and Rozansky noted, the Euler characteristic of this
complex is the HOMFLY polynomial of links (after a slightly odd change
of variables).  If $\bar\si$ is a knot, then the homology of this
complex is finite dimensional.

It is worth noting, we can extend the functor to singular braids in an
obvious way by defining a Rouquier complex $F(\si_i^!)$ for a single
intersection point between the $i$th and $i+1$st strand of a braid,
and applying the scheme above to a braidlike projection of the
singular knot.  This produces a categorification of the MOY state-sum
invariant of singular links.  This will not necessarily have
finite-dimensional homology, even for singular knots.

To obtain $\fr{sl}_N$-homology instead of HOMFLY homology, we need
only replace Hochschild homology by a slightly different functor.

The algebraic basis for the construction of this functor is the theory of matrix
factorizations.  Let $M$ be a $\Z$-graded module over a ring $S$.

\begin{defi}
A \emph{(}$\Z$-graded\emph{)} {\bf matrix factorization} on $M$ with
potential $\vp\in  S$ is a map $d=\dpl+\dm:M\to M$ with $\dpm$ of
graded degree $\pm 1$ such that $d^2=\vp$.
\end{defi}

Though this is not the usual definition of a matrix factorization (where
typically we only assume a $\Z/2$ grading), this richer structure is
also useful from the perspective of knot theory. 

While this may look like a daunting definition, we will only be
interested in essentially a single example of a matrix factorization.

Fix an integer $N$, and index $1\leq i\leq n$, and let $\vp_i=x_i\otimes 1-1\otimes x_i$ and $\psi_i=\frac{x_i^N\otimes 1-1\otimes x_i^N}{x_i\otimes 1-1\otimes x_i}$. Define the matrix factorization $Z_i$ over $S\otimes S$ to be rank 2, with one copy of $S\otimes S$ in degree 0 and one in degree 1 with and $\dpl$ and $\dm$ by
\begin{equation*}
  Z_{i}=\xymatrix{ S\otimes S \ar@/^/[r]^{\vp_i}& S\otimes S\ar@/^/[l]^{\psi_i}}.
\end{equation*}
and let $Z=\otimes_{i=1}^n Z_i$ (where the tensor product is essentially that of complexes applied to both differentials).  Note that the potential of $Z$ is $\sum_{i=1}^nx_i^{N+1}\otimes 1-1\otimes X_i^{N+1}$.

If $M$ is an $S-S$ bimodule annihilated by $p\otimes 1-1\otimes p$ for
any symmetric polynomial $p\in S^{S_n}$, then $Z\otimes_{S\otimes S}
M$ is matrix factorization of potential 0.  That is, the total
differential $d$ is an honest differential, so the total homology
$\tilde{\HH}_N(M)=H(Z\otimes_{S\otimes S}M,d)$ is well defined, and in fact
carries a single grading, which is a linear combination of the polynomial and matrix factorization gradings on $Z\otimes_{S\otimes S}M$ in which $d$ is homogeneous.

Note that if we only consider $\dm$, we obtain a free resolution of
$S$ as a bimodule over itself.  Thus the homology of the complex
$H^i(Z_-\otimes_{S\otimes S} M)$ is simply $\HH^i(M)$ for any $S-S$-bimodule $M$.  Thus, we can think of $\tilde{\HH}_N(M)$ as a sort of non-flat deformation of $\HH^*(M)$.

\begin{defi}
  The {\bf $\fr{sl}_N$ \KRH} $\KR_N(\bar{\si})$ of $\bar{\si}$ is the bigraded
  complex $\tilde{\HH}_N(F(\si))$.  This again has finite dimensional homology.
\end{defi}
In fact, Jacob Rasmussen has shown that the total dimension of this
homology if bounded above by that of the HOMFLY homology \cite{Ras06}.

The equivalence of this definition to that originally given in
\cite{KR1} was communicated to the second author by Mikhail Khovanov,
but to the best of our knowledge the first full proof appeared in
\cite{Web06a}.

\section{Categorified Vassiliev derivative.}
\label{sec:the categorification of the Vassiliev derivative}

In this section, we will discuss a general schema for Vassiliev theory of knot homology, as introduced in the introduction and \cite{Shi07a}.  

Let $\A$ be an abelian category, and let $\KA$ be the category of
complexes in $\A$ with morphisms considered up to homotopy.  The
category $\KA$ is not abelian; it no longer makes sense to consider
the kernel or cokernel of a map.  The closest notion we have is that
of the cone of a map.
 
As usual, for a complex $X=(X^i,d_X^i)$, define the shift $X[j]$ of
$X$ by $$(X[j])^i=X^{i+j}\hspace{10mm} d_{X[j]}=(-1)^jd_X$$
 
\begin{defi}
  Let $f:X\rightarrow Y$ be a chain morphism. The {\bf cone} of f is
  the complex
  \begin{equation*}
C_f^i=X[1]^i \oplus Y^i\hspace{10mm} d_{C_f}(x^{i+1},y^i)=(-d_X x^{i+1},f(x^{i+1})-d_Yy^i)
\end{equation*}
\end{defi}

In $\KA$, the cone naturally fits into an exact triangle:
  $$
\xymatrix{
&C_f\ar[dl]^w_{[1]}&\\
X^\bullet \ar[rr]^f&&Y^\bullet \ar[ul]^v
} \label{eq:tri1}
$$

For our purposes, cones have two important properties.  The first is a rather trivial observation.  If we let $\chi(X)$ be the Euler characteristic of $X$ (a class in the Grothendieck group $K^0(\A)$ then 
\begin{proposition}\label{add-Eu}
  $\chi(C_f)=\chi(X)-\chi(Y)$
\end{proposition}

The second is the behavior of successive cones.  Whenever we have a commuting square of chain maps
\begin{equation*}
  \xymatrix{X_{-,-}\ar[r]^{\vp_-}\ar[d]_{\psi_-}&X_{-,+}\ar[d]^{\psi_+}\\X_{+,-}\ar[r]_{\vp_+}&X_{-,-}}
\end{equation*}
then one has natural induced maps $\vp:C_{\psi_-}\to C_{\psi_+}$ and
$\psi:C_{\vp_-}\to C_{\vp_+}$.  It is a simple exercise to show that
$C_{\psi}$ and $C_\vp$ are naturally isomorphic.  More generally if we
have a commuting hypercube of any dimension, we will get the same
answer taking cones in any order.  This iterated cone can be seen
as the total complex of a bicomplex defined by our chain maps.

Consider a point of self intersection of the discriminant of
codimension $n$. There are $2^n$ chambers adjacent to this
point. Since the discriminant was resolved by Vassiliev \cite{Vas92}, this
point can be considered as a point of transversal self intersection of
n hyperplanes in $R^n$, or an origin of the coordinate system of
$R^n$.

For an invariant $\lambda$ of knots valued in an abelian group, we can
extend $\la$ to the discriminant by the Vassiliev derivative.  If $K$ is
a singular knot with $n$ self-intersection points, then there are $n$
codimension 1 ``walls'' of the discriminant intersecting transversely at
$n$, each being cooriented (having a ``positive'' and a ``negative''
side).  Thus, any neighborhood of $K$ is split into $2^n$ chambers, one
for each map $\si:[1,\ldots, n]\to\{\pm 1\}$.  Let $K_\si$ be a
representative for chamber corresponding to $\si$, and let
\begin{equation*}
\la(K)=\sum_{\si}(-1)^\si \la(K_\si)
\end{equation*}
where $\nu(\si)=\sum_{i=1}^n(1+\si(i))/2$.

As usual in categorification, when we pass to a categorification, we
would like to replace the sum above with a chain complex.  The
previous work of the first author \cite{Shi07a} suggests that to
categorify a knot homology theory $\HE$, we should define {\df
  wall-crossing morphisms} $\wc_w:\HE(K_-)\to \HE(K_+)$ where $K_\pm$
are the knots adjacent to a wall $w$, such that the diagram
\begin{equation}\label{wall-com}
\xymatrix{\HE(K_{-,-}) \ar[r]^{\wc_{w_1}}\ar[d]_{\wc_{w_2}}&\HE(K_{-,+})\ar[d]^{\wc_{w_2}}\\
\HE(K_{+,-})\ar[r]_{\wc_{w_1}}&\HE(K_{+,+})}
\end{equation}
commutes, where $K_{\pm,\pm}$ are the knots adjacent to a generic
element in the intersection of two walls $w_1$ and $w_2$.

There is a unique cobordism of genus 1 joining $K_-$ and $K_+$,
and one would hope to be able to define $\wc_w$ as simply the functor
$\HE$ applied to this cobordism.  Obviously, if $\HE$ is functorial on
the nose, this will work, but at the moment this approach is captive
to the sign problems which appear in many knot homology theories,
including Khovanov-Rozansky homology.  In this sense, one can consider
the existence of wall-crossing morphisms as a weaker version of
fixing the signs of functoriality.

Now, assume that we have constructed such morphisms $\wc_w$ for all
walls in the discriminant. At each singular knot $K$, we can define
$\HE(K)$ as an iterated cone $C_{\wc_{w_1},\ldots,\wc_{w_n}}$  of the wall-crossing maps corresponding to walls $w_1,\ldots,w_n$ containing $K$.

As an alternative description, we can construct a complex $\mc
C_K$ such that
\begin{equation*}
(\mc C_K)_i= \bigoplus_{\nu(\si)=i}\HE(K_\si)   
\end{equation*}
and differentials given by appropriate sums of wall-crossing maps (as
usual, we will need to add signs to make sure we have a complex, but
this can be done by the standard conventions of supermathematics), by
simply collapsing the grading on the hypercube with each resolution at
a corner.

This is now a double-complex in the category of matrix factorizations of
potential 0.  We can then extend $\HE$ to the discriminant by taking
the total complex
\begin{equation*}
\HE(K)=\mathrm{Tot}(\mc C_K).
\end{equation*}

Proposition~\ref{add-Eu} implies 
\begin{corollary}
  The Euler characteristic of $\HE(K)$ is the extension of the knot
  invariant $\chi(\HE)$ to the discriminant, that is
  \begin{equation*}
    \chi(\HE(K))=(-1)^{\nu(\si)}\chi(K_\si),
  \end{equation*}
\end{corollary}

\section{Wall-Crossing Morphisms}
\label{sec:wall-cross-morph}

While the above discussion covered an essentially formal situation which
could apply to any categorification, we still need to define the
wall-crossing maps themselves, which will require us getting our hands
(a little) dirty.

Fix a wall $w$ and a generic singular knot $K$ in $w$ and let $K_\pm$
be the knots on its positive and negative sides.  Fix a braid-like
projection of $K$.  In this projection, $K$ has a single
self-intersection point, and a projection of $K_+$ (resp. $K_-$) is
obtained by resolving this self-intersection point to a positive
(resp.  negative) crossing. Furthermore, $K$ is the closure of a
singular braid $\si_!$.  Let $\si_+,\si_-$ be the resolution of the
self-intersection point of $\si_!$.  We may assume (by moving the
cutting point) that for some braid $\beta$, we have $\si_!=\si_i^!\beta, \si_+=\si_i\beta$
and $\si_-=\si_i^{-1}\beta$.  This in turn implies that
\begin{equation*}
  F(\si_+)\cong F(\si_i)\otimes_S F(\beta)\hspace{10mm}F(\si_-)\cong F(\si_i^{-1})\otimes_S F(\beta)
\end{equation*}
\begin{theorem}
  There exists a natural (up to scalar) chain map
  $\wc:\KR_N(K_-)\to\KR_N(K_+)$ for $N=1,\ldots,\infty$.  If $K$ is an
  $n$-singular link, there is consistent choice of scalars so that the
  induced cube of wall-crossing maps is commutative, and thus the
  iterated cone on this cube is independent of these scalars.
\end{theorem}
\begin{proof}
  This map is induced by an element in $\Ext^1(F(\si_-),F(\si_+))$,
  which, in turn, comes from one in $\vp\in\Ext^1(F(\si_i^{-1}),
  F(\si_i))$.  This element can be described in several ways.

  If $\ti S$ is the quotient of $S\otimes S$ by the ideal generated by
  $x_j\otimes 1 -1\otimes x_j$ for $j\neq i,i+1$, then we have a
  natural exact sequence of complexes:
  \begin{equation*}
\xymatrix@C=3cm{\ti S\ar[r]^{x_i\otimes 1 -1\otimes x_i}&\ti S\{-2\} \ar[r]& S\{-2\}\\ \ti S\{2\}\ar[r]_{x_ix_{i+1}\otimes 1-1\otimes x_ix_{i+1}} \ar[u]_{x_{i}\otimes 1 -1\otimes x_{i+1}} &\ti S\{-2\}\ar[r] \ar[u]_1& S_i\{-1\}\ar[u]_{m_i}}
  \end{equation*}

  By simple degree considerations, this exact sequence induces an
  injection $\mathrm{Hom}(\{\ti S\{2\}\too \ti S\}, F(\si_i))\hookrightarrow
  \Ext^1(F(\si_i^{-1}),F(\si_i))$.  Since there is a unique (up to
  scalar) projection map $\{\ti S\{2\}\too \ti S\} \too
  F(\si_i)$, we can take $\vp$ to be the image of this.

   Alternatively, if we let $S'=\ti S/(x_i\otimes 1 -1\otimes x_i)^2(x_i\otimes 1 -1\otimes x_{i+1})$, then the exact sequence
  \begin{equation*}
    \xymatrix@C=3cm{S_i\{1\}\ar[r]&S'\{-2\} \ar[r]& S\{-2\}\\ S\{2\}\ar[r] \ar[u]_{\iota_i} &S'\{-2\}\ar[r] \ar[u]_{1}& S_i\{-1\}\ar[u]_{m_i}}
  \end{equation*}
  is a realizaion of the class $\vp\in\Ext^1(F(\si_i^{-1}),F(\si))$.

  By standard homological algebra, any element $\psi\in \Ext^1_R(M,N)$
  induces a canonical map $\HH^i(\psi):\HH^i(M)\to\HH^{i+1}(N)$.  Thus, the
  image of $\vp$ in $\Ext^1(F(\si_-),F(\si_+))$ induces a map
  $\wc:\KR_N(K_-)\to\KR_N(K_+)$, which is our wall-crossing map.

  Note that if $K$ is an $n$-singular link, we can define the
  wall-crossing element of $\Ext^1$ in a different tensor factor, so
  they will commute in the Yoneda product, and thus induce commuting
  maps on $\HH^i$.  Even though $\vp_\theta$ for each singular point
  $\theta$ is only defined up to a scalar, changing the $\Ext^1$-term
  in the factor corresponding to $\theta$ by a scalar will change all
  the wall-crossing maps for that wall by the same scalar, so the
  iterated cone is still well-defined.
\end{proof}

Let $K$ be a singular
link, 
and $\{K^{\la}\}$ be the collection of resolutions of $K$. We have a
cube of wall-crossing maps whose vertices are the complexes
$\KR_N(K^{\la})$ as $\la$ ranges over resolutions.  We denote the
total complex of multi-complex by $\VKR_N(K)$, with the grading
inherited from $\KR_N(K^\la)$ (or bigrading in the case of
$\KR_\infty(K)$).  The complex $\VKR_N(K)$ can also be realized as an
iterated cone, over the wall-crossing maps, for the various walls
which $K$ lies on.

\begin{defith}
  The homology of $\VKR_\infty(K)$ is a triply graded homology theory
  for singular links and its Euler characteristic is the Vassiliev
  derivative of the HOMFLY polynomial.

  Similarly, the homology of $\VKR_N(K)$ for $N<\infty$ is a
  categorification of the Vassiliev derivative of the
  $\fr{sl}_n$-quantum invariants.
\end{defith}
\begin{proof}
  The proof of invariance simply follows that of invariance of HOMFLY
  homology based on the Markov moves in \cite{Web06a}.  Markov I is
  clear, because $\vp$ is depends on a single one of the tensor
  factors which are cyclically permuted by Markov I.  For Markov II,
  we need only check that the inclusion of
  $\mathrm{Hom}(\EuScript{F}^j (\si_i^{-1}), \EuScript{F}^{j+1}(\si_i))$
  to $\mathrm{Hom}(\EuScript{F}^j (\si_-), \EuScript{F}^{j+1}(\si_+))$
  matches that after the stabilization under the isomorphisms that hold for all $\si\in B_n$.
  $\EuScript{F}^j(\si)\cong \EuScript{F}^j(\si\si_{n})$ and
  $\EuScript{F}^j(\si)\cong \EuScript{F}^{j-1}(\si\si_{n}^{-1})$.
  Both these isomorphisms are induced by spectral sequences for the
  tensor product $F(\si)\otimes F(\si_{n}^{\pm})$, and by the
  functoriality of these spectral sequences, the inclusions coincide.
\end{proof}

\begin{proposition}
  The homology of $\VKR(K)$ is finite dimensional for any singular knot $K$.
\end{proposition}
\begin{proof}
  The iterated cone is the total chain complex of a double complex
  whose horizontal components are the complexes $\KR_N(K^\la)$, whose
  homology is finite dimensional (see, for example, \cite[Propositon
  7.1]{Ras06}).  The spectral sequence of a double complex shows that
  the homology of the total complex is finite-dimensional as well.
\end{proof}

\def\cprime{$'$}
\providecommand{\bysame}{\leavevmode\hbox to3em{\hrulefill}\thinspace}
\providecommand{\MR}{\relax\ifhmode\unskip\space\fi MR }
\providecommand{\MRhref}[2]{%
  \href{http://www.ams.org/mathscinet-getitem?mr=#1}{#2}
}
\providecommand{\href}[2]{#2}

\end{document}